\documentclass[12pt]{amsart}
\usepackage{graphicx}
\usepackage{amsmath}
\usepackage{amstext}
\usepackage{amsfonts}
\usepackage{amsthm}
\usepackage{amssymb}

\newtheorem{theorem}{Theorem}
\newtheorem{definition}{Definition}

\newtheorem{lemma}{Lemma}

\newcommand{\Ar}{\mathbb{R}}
\newcommand{\Be}{\frak{B}}

\newcommand{\CSD}{\mathcal{L}}

\newcommand{\ep}{\epsilon}

\newcommand{\Bord}{\Omega F}

\newcommand{\HF}{HF}

\newcommand{\term}[1]{\textbf{#1}}

\newcommand{\Loopc}{\Lambda(T^*M)}
\newcommand{\Loop}{\Lambda(M)}
\newcommand{\Loopv}{\Lambda(TM)}

\begin{document}

\title{Semi-Infinite Cycles in Floer Theory: Viterbo's Theorem}
\author{Max Lipyanskiy}
\date{}
\address{Department of Mathematics, Columbia University, New York, NY 10027 }
\email{mlipyan@math.columbia.edu}
\maketitle

\begin{abstract} This is the first of a series of papers on foundations of Floer theory.  We give an axiomatic treatment of the geometric notion of a semi-infinite cycle.  Using this notion, we introduce a bordism version of Floer theory for the cotangent bundle of a compact manifold $M$. Our construction is geometric and does not require the compactness and gluing results traditionally used to setup Floer theory.  Finally, we prove a bordism version of Viterbo's theorem relating Floer bordism of the cotangent bundle to the ordinary bordism groups of the free loop space of $M$.   
\end{abstract}

\section{Introduction}
In the mid 1980's, Andreas Floer developed a new homology theory for infinite dimensional manifolds, based on a version of Morse theory.  As an immediate application, Floer was able to resolve a version of the Arnold conjecture on the existence of periodic orbits of a Hamiltonian vector field on a symplectic manifold \cite{Floer1}.  Floer's theory, however, was very different from the infinite dimensional Morse theory developed by Morse and Bott for energy functional on a Riemannian manifold.  For instance, the Hessian of the functional that Floer considered had infinitely many positive and negative eigenvalues.  Therefore, in his setting, passing a critical point was topologically equivalent to adding a handle along an infinite dimensional sphere.  Since such spheres are contractible, it was not all clear what geometric/topological properties of the underlying space this homology theory captured.    \\\\
To our knowledge, Atiyah was first to suggest that Floer's construction may be formulated as geometric theory of semi-infinite cycles for the given infinite dimensional manifolds \cite{Ati}.  Let us briefly recall his suggestion.  While Atiyah's article discussed the case of instanton Floer theory, we prefer to formulate it for the case of symplectic Floer homology.  Let $M$ be a closed smooth manifold and let $T^*M$ be its cotangent bundle.   As our configuration space we take  $\Loopc$, the space of free loops on $T^*M$.  By completing the loop space with respect to some Sobolev norm, we may view this loop space as a Hilbert manifold.  As we will explain in the paper, $\Loopc$ inherits a polarization  in the form of a splitting $$T\Loopc = T^+\Loopc\oplus T^-\Loopc$$  Atiyah's proposal is that Floer homology associated to $\Loopc$, may be defined by considering smooth mappings $$\sigma: P\rightarrow \Loopc$$ where $P$ is some Hilbert manifold hitting roughly half of the $T\Loopc$.  More precisely, $D\sigma$ should be Fredholm when projected to $T^-\Loopc$ and compact when projected to $T^+\Loopc$.  This is the origin of the term semi-infinite.  We can also reverse the parts of the polarization and consider maps $$\tau:Q\rightarrow \Loopc$$  which are Fredholm when projected to $T^+\Loopc$ and compact when projected to $T^-\Loopc$.  Thus, on an intuitive level, we arrive at a notion of ``upper'' and ``lower'' semi-infinite cycles.\\\\
Now, the key point is that there is a paring between these two types of cycles.   Assuming that $\sigma$ is transverse to $\tau$, we can form the fibre product $$\sigma \cap \tau:P\cap Q \rightarrow \Loopc$$  Here $P\cap Q$ is the set of points $(p,q)\in P\times Q$ with $\sigma(p)=\tau(q)$.  A simple argument from functional analysis  shows that $P\cap Q$ must be finite dimensional.  Thus, $\sigma \cap \tau$ defines a ordinary cycle in $\Loopc$.  Therefore, one may hope to detect nontrivial semi-finite cycles by showing that their pairing gives rise to nonzero ordinary homology classes in $\Loopc$.  \\\\
The crucial point missing from this discussion is the issue compactness of $P\cap Q$.  Since $P$ and $Q$ can never be compact, it is not clear how to restrict the class of allowed maps $\sigma$ and $\tau$ to ensure that $\sigma\cap \tau$ is compact.  For instance, given $\sigma$, we need to rule out $$\Sigma: P\times [0,1) \rightarrow \Loopc$$  with $\Sigma(p,t)=\sigma(p)$.  Otherwise, every cycle will be a boundary.\\\\
In this work, we propose a set of axioms for $\sigma$ and $\tau$ as above, designed to ensure compactness of $\tau \cap \sigma$.  Analogous to the work of Floer, a necessary ingredient is the symplectic action functional $$\CSD:\Loopc \rightarrow \Ar$$ We will require $\sigma$ and $\tau$ to satisfy certain point-set topological constraints with respect to $\CSD$.  As we shall demonstrate, such constrains will ensure that the fibre products are always compact.  We defer the precise definitions to the main body of the paper.  We would like to point out that,  unlike Floer's construction, we do not need to assume that $\CSD$ has nondegenerate critical points.  In fact, $\CSD$ need not be differentiable for the groups to be defined.  To summarize, the key point in constructing a nontrivial semi-infinite homology theory based on geometric cycles is to develop some notion of semi-infinite topology.\\\\
Given an appropriate notion of semi-infinite cycles, the simplest groups to construct are semi-infinite versions of bordism groups $\Bord_*(\Loop)$ and $\Bord^*(\Loop)$.  Taking intersections, we will obtain a pairing: $$\Bord_*(\Loop)\otimes \Bord^*(\Loop)\rightarrow \Omega_*(\Loop)$$  Here $\Omega_*(\Loop)$ denotes ordinary oriented bordism.  \\\\
Given the existence of Floer bordism, it is natural to try to express them in terms of more familiar invariants.  Let us recall the following result, originally due to Viterbo (see \cite{Viterbo1}, \cite{SalWeb} and \cite{AbbonSchwa}):
\begin{theorem}
There exists an isomorphism $\HF_*(T^*(M))\cong H_*(\Loop)$.
\end{theorem}
Now, we state the main result of this paper which is the bordism analogue of Viterbo's theorem.  Let $$\pi:\Loopc \rightarrow \Loop$$ be the projection to the base.   If $$\sigma:P\rightarrow \Loop$$ is an element of $\Omega_*(\Loop)$, we can form the pullback $$\pi^*(\sigma):\pi^{-1}(P)\rightarrow \Loopc$$  
\begin{theorem}
\label{Main}
$\pi^*$ induces an isomorphism $\Omega_*(\Loop) \cong \Bord_*(\Loopc) $. 
\end{theorem}   

The proof of this theorem is rather different from the existing proofs of the result of Viterbo since we do not make any use of Morse theory.  In fact, as mentioned above, the functional defining Floer bordism need not have nondegenerate critical points.  As we shall demonstrate, the proof boils down to a  homotopy argument where the work consists of checking that the homotopies satisfy all the topological assumptions imposed on the cycles. As such, it completely avoids the use of partial differential equations. \\\\
Finally, let us point out that since Atiyah's initial insight, Floer homology has been recast in terms of classical algebraic topology by several authors.  For instance, the work of Cohen, Jones and Segal in \cite{CohenJonesSegal} as well as the work of Cohen in \cite{CohenViterbo} uses the Morse-Floer gradient flow data to construct a spectrum associated to the configuration spaces. In Seiberg-Witten theory, a formulation of Floer theory in terms of spectra has been carried out by Manolescu in \cite{Manolescu1}.  \\\\
The main ideas of this work stem from the author's thesis written under the supervision of Tomasz Mrowka.  The idea of developing axioms for semi-infinite chains has its roots in unpublished work of Tomasz Mrowka and Peter Ozsvath from the late '90.  I am happy to use this opportunity to thank Tomasz for his sharing his thoughts on this subject.  In addition, I would like to thank Peter Kronheimer, Cliff Taubes as well as Dennis Sullivan for several very helpful conversations as well as their interest in this project.        
        
\section{Bordism Groups of the Loop Space}
In this section we briefly recall the construction of the Hilbert manifold structure on $\Loop$ and set up some notation.  For the entirety of this paper, let $M$ be a closed, smooth  manifold of dimension $m<\infty$.  For convenience, we fix an embedding $M\subset \Ar^N$.  Also, we will make use of a Riemannian metric on $M$ which we will denote by $g$.  We let $\nabla$ be the associated Levi-Civita connection on $T(M)$.  Given a loop $$\gamma:S^1\rightarrow M$$   we have the energy, $$E(\gamma)=\frac{1}{2}\int_{S^1}|\gamma'|^2_g$$
Given an smooth loop $\gamma$, let us define the square of the $L^2_1$-norm as $$\int_{S^1}|\gamma|^2+E(\gamma)$$  Here, $|\gamma|$ is defined using the norm on $\Ar^N$.  
\begin{definition}
Let $\Loop$ denote the $L^2_1$-completion of the set of smooth loops $\gamma$.
\end{definition}
As is well known (see \cite{Floer1} for example), $\Loop$ is actually a Hilbert manifold.  We may describe the smooth structure on $\Loop$ as follows.  Given a smooth $\gamma$, let $\Gamma(\gamma^*(TM))$ be the space of $L^2_1$ sections of $\gamma^*TM$. $\nabla$ on $TM$ gives rise to a Hilbert space structure on $\Gamma(\gamma^*(TM))$ as follows.  Given $\eta\in \Gamma(\gamma^*(TM))$, let $$||\eta||_{L^2_1}^2=\int_{S^1} |\eta|^2_g+ \int_{S^1} |\nabla_{\partial_t} \eta|^2_g$$  We let $exp_\gamma(\eta)(t)\in M$ be point of $M$ given by the time one map of the  unique geodesic starting at $\gamma(t)$ with initial speed $\eta(t)$.  We have:
\begin{lemma}
Let $\ep>0$ and take $\Gamma(\gamma^*TM)_\ep$ to be sections with $||\eta||_{L^2_1} < \ep$.  For small enough $\ep$, $$exp_\gamma:\Gamma(\gamma^*(TM))\rightarrow \Loop$$ is a homeomorphism to its image.  Furthermore, as we vary $\gamma$, such charts provide $\Loop$ with the structure of a smooth Hilbert manifold.
\end{lemma}   
We now recall the construction of the bordism group of $\Loop$.  Consider smooth maps $$\sigma:P\rightarrow \Loop$$ where $P$ is an oriented, closed, smooth finite dimensional manifold.  Two such maps $\sigma_1:P_1\rightarrow \Loop$, $\sigma_2:P_2\rightarrow \Loop$ are isomorphic if there exists a orientation preserving isomorphism $f:P_1\rightarrow P_2$ such that $\sigma_1=\sigma_2\circ f$. We denote the isomorphism by $\cong$. An element, $\sigma$  is \textbf{trivial} if there exists an oriented, finite dimensional manifold with boundary $Q$, a smooth map $$\tau:Q\rightarrow \Loop$$ and a orientation preserving diffeomorphism $$\phi:\partial P\rightarrow Q$$ such that $\sigma=\tau \circ \phi$.  
\begin{definition}
$\sigma_1 \sim \sigma_2$ if there exists trivial cycles $\sigma'_1$ and $\sigma'_2$ such that $\sigma_1\sqcup \sigma_1'\cong \sigma_2\sqcup \sigma_2'$.   
\end{definition}     
We have the following elementary lemma:
\begin{lemma}
If $\sim$ is an equivalence relation.
\end{lemma}
\begin{proof}
The only part that needs to be checked is the transitive property.  If $\sigma_1\sqcup \sigma_1' \cong \sigma_2 \sqcup \sigma_2'$ and $\sigma_2\sqcup \sigma_2'' \cong \sigma_3 \sqcup \sigma_3''$ where $\sigma'_i$, $\sigma''_i$ are trivial, we get $$\sigma_1\sqcup \sigma_1' \sqcup \sigma''_2 \cong \sigma_2 \sqcup \sigma_2'\sqcup \sigma_2'' \cong \sigma_3 \sqcup \sigma_3'' \sqcup \sigma_2'$$ 
\end{proof}
\textbf{Remark.}  The tubuluar neighborhood theorem implies that if $\sigma$ is equivalent to a trivial cycle it is itself trivial.
\begin{definition}
Let $\Omega_*(\Loop)$ be the group generated by maps $\sigma:P\rightarrow \Loop$ modulo the equivalence relation.  The additive structure is given by disjoint union. The grading is given by the dimension of $P$.
\end{definition}

There is a geometric construction of cobordism due to Daniel Quillen in \cite{Quillen1} (also see \cite{BakerOzel} for the infinite dimensional case).  For this, consider smooth Fredholm maps $\tau :Q\rightarrow \Loop$, where $Q$ is a separable Hilbert manifold.  Recall (see \cite{KM} for an exposition) that there is a determinant line bundle $det(\tau)\rightarrow Q$  with the fiber over $q$ given by $$\Lambda^{max} ker(D_q\sigma)\otimes (\Lambda^{max}coker(D_q\sigma))^*$$ 
\begin{definition}
An \textbf{orientation} of $\tau$ is a choice of orientation of $det(\tau)$. 
\end{definition}
If $Q$ is a manifold with boundary, an orientation of $\tau$ induces an orientation of $\tau_{|\partial Q}$ using the "outward normal first" convention.  We define the cobordism group as follows.  We consider smooth, proper, oriented Fredholm maps  $$\tau:Q\rightarrow \Loop$$ where $Q$ is a Hilbert manifold without boundary.  We assume that each component of $\tau$ has the same Fredholm index.  We call such a map \textbf{trivial} if there exists a smooth, proper, oriented, Fredholm map $$\tau':Q'\rightarrow \Loop$$ where $Q'$ is a Hilbert manifold with boundary and a diffeomorphism $$\phi:Q\rightarrow \partial Q'$$ with $\tau=\tau' \circ \phi$ and $det(\tau)\cong det (\tau' \circ \phi)$ as oriented bundles.  Just like bordism, this leads to an equivalence relation on cycles: 
\begin{definition}
$\tau_1 \sim \tau_2$ if there exists trivial cycles $\tau'_1$ and $\tau'_2$ such that $\tau_1\sqcup \tau_1'\cong \tau_2\sqcup \tau_2'$.   
\end{definition} 
\begin{definition}
Let $\Omega^*(\Loop)$ be the group generated by such $\tau$ modulo $\sim$.  The group is graded by the Fredholm index of $\tau$.  
\end{definition}   
Finally, let us note that standard transverality arguments give rise to a pairing $$\Omega_a(\Loop)\otimes \Omega^b(\Loop)\rightarrow \Omega_{a+b}(\Loop)$$
The compactness of the intersection follows from the fact that elements of $\Omega^*(\Loop)$ are proper.

\section{Axioms for Floer Bordism}
\subsection{Basic Construction}
In this section we lay down the axioms necessary to set up a geometric version of Floer theory.  While this paper deals primerely with the example of the cotangent bundle of a closed manifold, other examples will be discussed in the sequel.  As described in the introduction, our main objective is to ensure impose certain point-set topological restrictions to ensure compactness of intersections.  Let $\Be$ be a separable Hilbert manifold.  The example to keep in mind is $\Be=\Loopc$. We let $T(\Be)$ denote the tangent bundle with its induced Hilbert space structure.  Note that at this point we do not specify an inner product on the tangent bundle.\\\\
We will need the notion of polarization of a Hilbert manifold. See \cite{CohenJonesSegal} and \cite{Douglas} for a general discussion as well as a topological classification.  For any open set $U\subset \Be$, let $T(U)$ be the restriction of $T(\Be)$ to $U$.        
\begin{definition}
A \textbf{polarizing chart} is a bundle isomorphism $$\phi: T^+(U)\oplus T^-(U)\cong T(U)$$ where $T^{\pm}(U)$ is a closed infinite dimensional subbundle of $T(U)$.  Let $\pi^{\pm}$ denote the projections $$\pi^{\pm}:T(U)\rightarrow T^{\pm}(U)$$ induced by the chart.  Two charts $\phi_i$, $\phi_j$ are \textbf{compatible} if $$\pi^{\pm}\circ \phi_j^{-1} \phi_i: T^{\pm}_i(U_i\cap U_j)\rightarrow T^{\pm}_j(U_i\cap U_j)$$ is Fredholm and $$\pi^{\mp}\circ \phi_j^{-1} \phi_i: T^{\pm}_i(U_i\cap U_j)\rightarrow T^{\mp}_j(U_i\cap U_j)$$ is compact.  
\end{definition}
\begin{definition}
A \textbf{polarization} of $\Be$ is a maximal compatible atlas. 
\end{definition}
\begin{definition}
A polarizing atlas is said to be \textbf{absolute} if for any two charts $\phi_i$, $\phi_j$, $$\pi^{+}\circ \phi_j^{-1} \phi_i: T^{\pm}_i(U_i\cap U_j)\rightarrow T^{\pm}_j(U_i\cap U_j)$$ has Fredholm index 0.  An \term{absolute polarization} is a maximal compatible absolute atlas. 
\end{definition}
We are ready to define the main object of this paper:
\begin{definition}
A \textbf{Floer Space} $(\Be,\CSD)$ consists of the following data:\\\\
1.  A Hilbert manifold $\Be$ together with a choice of a weak topology $T'$.  $T'$ is assumed to be more coarse then the ordinary (strong) topology of $\Be$.\\
2.  A polarization of $\Be$.\\
3.  A functional $\CSD:\Be \rightarrow \Ar$ which is continuous with respect to the ordinary topology.
\end{definition}
\begin{definition}
Given a Floer space $(\Be,\CSD)$, let $(-\Be, -\CSD)$ be the Floer space obtained by changing the sign of $\CSD$ and reversing the roles of $T^+(\Be)$ and $T^-(\Be)$.
\end{definition}
Let $P$ be a separable Hilbert manifold (possibly with boundary) and let  $\sigma:P\rightarrow \Be$ be a smooth map.  Our goal is to define the notion of semi-infinite cycle.  To motivate the definition, recall the following basic facts about Hilbert spaces:
\begin{definition}
A sequence $v_i$ in a Hilbert $\mathbb{H}$ is weakly convergent if for any $w\in \mathbb{H}$ we have:
 $$\lim_{i\rightarrow \infty} \langle w,v_i \rangle =\langle w,v_\infty \rangle$$
\end{definition}
\begin{lemma} \label{complemma}
A sequence $v_i$ is weakly precompact iff $|v_i|$ is bounded.  If $v_i$ is weakly convergent, $|v_\infty|\leq \liminf |v_i|$.  Finally, if $\lim |v_i|=|v_\infty|$, $v_i$ converge strongly to $v_\infty$.    
\end{lemma} 
Our definition of a semi-infinite cycle is basically the statement that $\CSD$ behaves like the negative of a Hilbert space norm when restricted to the cycle:
\begin{definition}
We say that $\sigma:P\rightarrow \Be$ is a \textbf{semi-infinite cycle} if the following axioms are satisfied:\\\\
Axiom 1.  $\CSD$ is bounded above on the image of $\sigma$.  If $y_i\in im(\sigma)$ converges weakly to $z\in \Be$, we have $$\liminf(\CSD(y_i))\leq \CSD(z)$$  \\
Axiom 2. If $y_i \in im(\sigma)$ has $\CSD(y_i)$ uniformly bounded, $y_i$ has a weakly convergent subsequence in $\Be$.\\\\
Axiom 3.  If $y_i=\sigma(p_i)$ converges weakly to $z\in \Be$ and $$\lim \CSD(y_i)=\CSD(z)$$ then $p_i$ has a convergent subsequence in $P$.\\\\
Axiom 4.  Let $\phi_i$ be a compatible polarizing chart for $\Be$. For each $p\in P$, we have $$\pi^-  \circ D\sigma_p:T(P)\rightarrow T^-(\Be)$$ is Fredholm and $$\pi^+  \circ D\sigma_p:T(P)\rightarrow T^+(\Be)$$ is compact.
\end{definition}
An immediate consequence of the axioms is the  compactness of intersections.  
\begin{definition}
Given cycle $\sigma:P\rightarrow \Be$ for $(\Be, \CSD)$ and a cycle $\tau:Q \rightarrow \Be$ for $(-\Be, -\CSD)$ we let $$\sigma\cap \tau:P\cap Q \rightarrow \Be$$ be the map from $P\cap Q\subset (p,q)\in P\times Q$ with $\sigma(p)=\tau(q)$.  Note that $P\cap Q$ depends on the maps and not only on the manifolds $P$, $Q$.
\end{definition}
\begin{lemma}
Given $\sigma:P\rightarrow \Be$ a semi-infinite cycle for $(\Be,\CSD)$ and $\tau:Q\rightarrow \Be$ a semi-infinite cycle for $(-\Be,-\CSD)$, we have that $P\cap Q$ is compact.  
\end{lemma}
\begin{proof}
Axiom 1 implies that $\CSD$ is bounded above and below on the image of $P\cap Q$.  Axiom 2 implies that the image of $P\cap Q$ is weakly precompact.  Now, consider a weakly convergent sequence $y_i=\sigma(p_i)=\tau(q_i)$ in image of $P\cap Q$ with limit $z$.  By passing to a subsequence, we can assume that $\CSD(y_i)$ converges.  Axiom 1 implies that $\lim(\CSD(y_i))\leq \CSD(z)$ and $\lim(-\CSD(y_i))\leq -\CSD(z)$.  Therefore, $\lim \CSD(y_i)=\CSD(z)$.  Axiom 4 implies that both $p_i$ has a convergent subsequence.  By passing to this subsequence we can repeat the argument to conclude that $q_i$ has a convergent subsequence as well.  Therefore, any sequence in $P \cap Q$ has a convergent subsequence.  
\end{proof}
Let us verify that $P\cap Q$ is a finite dimensional manifold when $\sigma$ is transverse to $\tau$:
\begin{lemma}\label{dimlemma}
Suppose $\sigma \pitchfork \tau$.  We have that  $P\cap Q$ is a compact, finite dimensional smooth manifold.  Furthermore,  given $(p,q) \in P\cap Q$ and a polarizing chart $\phi$, the local dimension of $P\cap Q$ near $(p,q)$ is $ind(\pi^- \circ D\sigma_p)+ind(\pi^+ \circ D\tau_q)$. 
\end{lemma}
\begin{proof}
This is an application of the implicit function theorem for Hilbert spaces.  Locally, we have $\sigma:V\rightarrow U$ and $\tau: W\rightarrow U$ where $U, V$, and $W$ are open balls in a Hilbert space.  We may locally express $P\cap Q$ as $(\sigma-\tau)^{-1}(0)$.  In terms of the local polarization, $T^+_{\sigma(p)}\oplus T^-_{\sigma(p)}$ we have $$D(\sigma-\tau)=\pi^-\circ D\sigma -\pi^-\circ D\tau +\pi^+\circ D\sigma -\pi^+\circ D\tau$$  with $\pi^+\circ D\sigma -\pi^-\circ D\tau$ compact. 
Therefore, $D(\sigma-\tau)$ is Fredholm with $$ind(D(\sigma-\tau))=ind(\pi^-\circ D\sigma)+ind(\pi^+\circ D\tau)$$  Since $D(\sigma-\tau)$ is assumed to be surjective, this index is the dimension of $P\cap Q$.
\end{proof}
We conclude this section by giving two examples of Floer spaces.\\\\
\textbf{Example 1.} Let $\mathbb{H}^{\pm}$ be a separable Hilbert space.  Let $\Be=\mathbb{H}^+\oplus \mathbb{H}^-$ and $$\CSD(h^+,h^-)=|h^+|^2-|h^-|^2$$   Let the polarization $T^{\pm}(\Be)$ be defined by the subspaces $\mathbb{H}^{\pm}$.  Let the weak topology be the usual weak topology associated to the Hilbert space.   As an example of a cycle for $(\Be, \CSD)$ take $H^-\subset \Be$.  Similarly, an example of a cycle for $(-\Be,-\CSD)$ is given by $H^+\subset \Be$.  We have $H^+\cap H^-=\{ 0\}$.    \\\\
\textbf{Example 2.}  As a model of $S^{\infty}$ take $\{a_n \}$ with $a_n \in \mathbb{C}$ such that $\sum_n |n||a_n|^2< \infty$ and $\sum_n|a_n|^2=1$.  We have a free $S^1$ action on $S^\infty$ given by $e^{i\theta}\{ a_n\}= \{ e^{i\theta} a_n\}$.  Let $\Be=S^\infty/S^1$ and $$\CSD(\{ a_n\})=\sum_n n|a_n|^2$$  Note that $\CSD$ is $S^1$ invariant and hence descends to $\Be$.  Since $S^\infty$ is a subset of the Hilbert space $L^2_{1/2}$, it inherits a weak topology. It is closed under this topology since the weak $L^2_{1/2}$ convergence of $\{a_n \}$ implies strong $L^2$ convergence of the same sequence.  On $\Be$, we take the quotient weak topology. Let us denote elements of $\Be$ by $[\{ b_n\}]$.  We have a collection of charts for $\Be$ indexed by $\mathbb{Z}$.  Namely, take $U_k=\{ b_n\}_{n\neq k}$ with $$\sum_{n}|n||b_n|^2<\infty$$  We have the embedding $$\phi_k:U_k\rightarrow \Be$$ taking $ \{ b_n\}_{n \neq k}$ to $[\frac{\{ b_n \}}{\sqrt{\sum_n|b_n|^2}}]$ with $b_k=1$. \\\\
The polarization in each chart is given by $T^-U_k=\{b_n\}_{n\leq 0}$ and $ T^+U_k=\{b_n\}_{n> 0}$.  The transition functions are given by $$\phi_l^{-1}\circ \phi_k(\{ b_i\})= \{ c_i\}$$ where $c_i=b_i/b_l$ where $i\neq l$ and we set $b_k=1$.  One can check that this respects the polarization.  In this example, let our cycle $P\subset \Be$ be $$P=\{ [\{a_n \}]|a_n=0, \forall n>0 \}$$ Notice that $\CSD_{P}$ is essentially the function induced by the negative of the $L^2_{1/2}$ norm on $S^\infty$.  Similarly, as a cycle in $(-\Be,-\CSD)$, take $Q$ given by setting $a_n=0$ for $n< 0$.  We have $P\cap Q=\{pt \}$.  Note that this implies that the corresponding elements of the Floer bordism groups are nontrivial, given that we can perturb the cycles to be transverse.

\section{Floer Bordism For the (Co)Tangent Bundle}
\subsection{Verification of the Axioms}
We now discuss our main example of a Floer space.  Let $\pi:TM\rightarrow M$ be the projection.    We will need to consider a particular Sobolev completion of the space of loops on $TM$.  
\begin{definition}
Let $\Loopv$ be the vector bundle over $\Loop$ with the fiber over $\gamma$, the space of $L^2$ sections of $\gamma^*(TM)$.  We will denote the elements of $\Loopv$ by pairs $(\gamma(t), v(t))$ with $\gamma\in \Loop$ and $v(t)\in \Gamma(\gamma^*(TM))$.  
\end{definition}
The proof that $\Loopv$ is a Hilbert bundle over $\Loop$ is standard (see \cite{Floer1} for example).  We will outline some the basic steps in the proof to set up some notation.  Let $\gamma\in \Loop$ be a smooth loop.  Recall from section 2 that a chart for $\Loop$ is given by $\Gamma^\ep(\gamma^*(TM))$.  Given $\mu(t)\in \Gamma_{L^2}(\gamma^*(TM))$, and $\eta(t)\in \Gamma^\ep(\gamma^*(TM))$, let $P_{\eta}(v(t))$ be the parallel transport with respect to $\nabla$ of $v(t)$ along ${exp_\gamma(\tau\eta(t))}_{0\leq \tau \leq 1}$.  This gives a chart $$\Gamma^{\ep}(\gamma^*(TM))\oplus \Gamma_{L^2}(\gamma^*(TM))\rightarrow \Loopv$$ 
The connection $\nabla$ defines a splitting $$T_{v,x}(TM)=T_x^vM\oplus T_x^hM$$ where $T_x^vM$ denoted the fibre direction and $T_x^hM$ denoted the horizontal direction.  Similarly, the charts induced by $\nabla$ on $\Loopv$ give rise to a splitting of $$T\Loopv=T^+\Loopv \oplus T^- \Loopv$$  Here, $T^+\Loopv$ corresponds to $\Gamma^\ep(\gamma^*(TM))$ while $T^-\Loopv$ corresponds to the fiber $\Gamma_{L^2}(\gamma^*(TM))$.  Hence the connection $\nabla$ gives rise to an absolute polarization of $\Loopv$.  Note that in this case the polarization is defined by a global splitting of $T\Loopv$. \\\\
Aside from the usual manifold topology, $\Loopv$ has a weak topology induced by the weak $L^2_1$ topology on the base and the weak $L^2$ topology on the fiber.  More precisely, since $$i:M\rightarrow  \Ar^n$$ we have $$Ti: TM\subset \Ar^N \times T\Ar^N$$  We take the weak topology on $\Loopv$ to be the weak $L^2_1$ topology on the first factor and the weak $L^2$ topology on the second factor.  We may characterize this topology as follows.  Given a sequence $(\gamma_i,v_i)$ with a bound on the $L^2_1$ norm of $\gamma$ and the $L^2$ norm of $v$, we may pass to a subsequence with $\gamma_i$ converging to $\gamma_\infty$ in $C^0$.  Therefore, we may use the connection to identify all the fibers over $\gamma_i$ for $i$ large.  This way, we obtain a $L^2$ bounded sequence $v'_i$ over $\gamma_\infty$.  By the usual Hilbert space convergence, a subsequence of $v'_i$ converges weakly in $L^2$ to $v'_\infty$.  Since $\gamma_i$ converges to $\gamma_\infty$ in $C^0$, we have that $v_\infty'$ is the weak limit of $v_i$.     
\begin{definition}
Let $\CSD(\gamma,v)=\frac{1}{2}\int_0^1| \gamma'|^2-\frac{1}{2}\int_0^1|v|^2$
\end{definition}
We have specified all the necessary ingredients to define a Floer space  $(\Loopv, \CSD)$.  Since the polarization is given by a global splitting, in this case, we can refine the definition of bordism to take into account orientations and grading.  
\begin{definition}
A cycle $\sigma: P\rightarrow \Loopv$ is said to have index $n$ if $\pi^- \circ D_p \sigma$ has  index $n$ for all $p\in P$.  
\end{definition} 
 To discuss orientations, we recall the notion of the determinant bundle:
\begin{definition}
Given a cycle $\sigma$, let $det^-(\sigma)\rightarrow P$ denote the real line bundle with fiber $$\Lambda^{max} ker ( \pi^-\circ D_p\sigma) \otimes (\Lambda^{max} coker (\pi^-\circ D_p\sigma))^*$$ over $p\in P$.             
\end{definition}
See \cite{KM} for a general discussion of determinant bundles.
\begin{definition}
A cycle $\sigma: P\rightarrow \Loopv$ is said to be oriented if $det^-(\sigma)$ is oriented.  
\end{definition}
 For a general Floer space, an absolute index grading as well as orientations exists under suitable assumptions on the polarization.  Since the main example of this paper has a natural choice of global splitting, we will not present the general construction here (however, see \cite{MorseHomology}). \begin{definition}
 A cycle $\sigma:P\rightarrow \Loopv$ without boundary is said to be \term{trivial} if there exists a cycle with boundary $\sigma':P'\rightarrow \Loopv$ and a diffeomorphism $$\phi:P\rightarrow \partial P'$$ with $\sigma=\sigma' \circ \phi$ and $det^-(\sigma)\cong det^- (\tau \circ \phi)$ as oriented bundles.
\end{definition}
As in the case of the loop space, we can declare $\sigma_1 \sim \sigma_2$ when there exists trivial $\sigma'_1$ and $\sigma'_2$ with $\sigma_1 \sqcup \sigma_1' \cong  \sigma_2 \sqcup \sigma_2'$.    
\begin{definition}
Let $\Bord_n(\Loopv)$ (resp. $ \Bord^n(\Loopv)$) denote the group of oriented boundaryless cycles in $(\Loopv,\CSD)$ (resp. $(-\Loopv,-\CSD)$) of index $n$ modulo $\sim$.  We let $$\Bord_*(\Loopv)=\oplus_{n\in \mathbb{Z}}\Bord_n(\Loopv)$$ and $\Bord^*(\Loopv)=\oplus_{n\in \mathbb{Z}}\Bord^n(\Loopv)$. 
\end{definition}
We give an example of a cycle in $\Bord_0(\Loopv)$.  Given any $\gamma \in \Loop$, let $V_\gamma$ be the fiber over $\gamma$ in $\Loopv$.  We have that $$\CSD_{|V_\gamma}=-\frac{1}{2}|\cdot|_{L^2}^2$$  Therefore, $\CSD$ is just proportional to the negative of the norm on $V_\gamma$.  It follows that the verification of the axioms of a cycle are reduced to the compactness lemma $\ref{complemma}$. \\\\
We  also give an example of a cycle in  $\Bord^0(\Loopv)$. For this, take $$i_0: \Loop\subset \Loopv$$ embedded as the zero section.  This time, $$\CSD_{|\Loop}(\gamma)=E(\gamma)$$  Since $M$ is assumed compact, the energy is bounds the $L^2_1$ norm of $\gamma$.  Therefore, the verification of the axioms again reduces to the compactness lemma $\ref{complemma}$.
\begin{lemma}
Given $(P,\sigma) \in \Bord_a(\Loopv)$ and $(Q,\tau) \in \Bord_b(\Loopv)$ such that $\tau \pitchfork \sigma$ we have $$\sigma \cap \tau \in \Omega_{a+b}(\Loop)$$ 
\end{lemma}
\begin{proof}
$\sigma \cap \tau $ is compact by the axioms.  It has dimension $a+b$ by lemma $\ref{dimlemma}$.  We need to verify that $P\cap Q$ inherits an orientation.  Over, $P\cap Q$ we have the oriented line bundle $det(\sigma)\otimes det(\tau)_{|P\cap Q}$. It is induced by restricting $\pi^-\circ D\sigma_p \oplus \pi^+ \circ D\tau_q$ to $P\cap Q$.  This family of Fredholm operators is homotopic to the family    $D\sigma_p +  D\tau_q$ over $P\cap Q$.  In view of the transversality of $\sigma$ and $\tau$,   $D\sigma_p +  D\tau_q$ has no cokernel and the kernel is isomoprhic to the tangent space of $P\cap Q$.  Therefore, an orientation of $det(\sigma)\otimes det(\tau)_{|P\cap Q}$ gives rise to an orientation of $P\cap Q$.  
\end{proof}
\subsection{Legendre Transform}
Since we have fixed a metric on $M$, we have an induced metric on $T^*M$ as well as an isomorphism $T^*(M) \cong T(M)$. We will identify vectors with 1-forms using the metric in what follows.
\begin{definition}
Let $\Loopc$ be the Hilbert vector bundle over $\Loop$ whose fiber over $\gamma\in \Loop$ is the space of $L^2$ sections of $\gamma^*T^*M$.
\end{definition}
In the context of symplectic Floer theory, the more natural space to consider is the free loop space of the cotangent bundle, $\Loopc$.  In this section we define an diffeomorphism $L:\Loopv \rightarrow \Loopc$ that will induce a Floer space structure on $\Loopc$.
\begin{definition}
Let $L:\Loopv\rightarrow \Loopc$ be the diffeomorphism taking $(\gamma,v)$ to $(\gamma, v+\gamma')$.  
\end{definition}  
We have $$L_*(\CSD)(\gamma, w)=\CSD(\gamma, w-\gamma')=\frac{1}{2}\int_{S^1}|\gamma'|^2-\frac{1}{2}\int_{S^1}|w-\gamma'|^2=\int_{S^1}\langle w,\gamma' \rangle -\frac{1}{2}\int_{S^1}|w|^2$$  Since the isomorphism $L$ induces a polarization on $\Loopc$, we obtain a Floer space $$(\Loopc,L_*(\CSD))$$ It follows that $L$ induces an isomorphism of the corresponding bordism groups:
\begin{theorem}
We have isomorphisms $$\Bord_*(\Loopv)\cong \Bord_*(\Loopc)$$ 
 $$\Bord^*(\Loopv)\cong \Bord^*(\Loopc)$$ induced by $L$.
\end{theorem}

\subsection{Transversality}
In this section we will demonstrate that cycles can be perturbed to be transverse to the zero section $\Loop\subset \Loopv$.  The main idea is to apply Sard's theorem to a suitable family of perturbations.  We now turn to constructing a local model for our perturbations. Let $$\pi:\Loopv \rightarrow \Loop$$ be the projection.  As a local chart of $\Loopv$, take $\Gamma^{\ep}(\gamma^*_0(TM))\oplus \Gamma_{L^2}(\gamma^*_0(TM))$. This is a vector bundle over  $\Gamma^{\ep}(\gamma^*_0(TM))$.  We now define a section of $\Gamma^{\ep}(\gamma^*_0(TM))$ supported in $\Gamma^{\ep/2}(\gamma^*_0(TM))$.  Let $\rho:\Ar\rightarrow \Ar$ be a bump function equal to 1 near $0$ and vanishing outside $\ep^2/4$.  If $\mu(t)\in \Gamma_{L^2}(\gamma^*_0(TM))$, we may define our section as $$s_{\mu}(\eta)=\rho(||\eta||^2_{L^2_1})\mu$$  Here $\eta \in \Gamma^{\ep}(\gamma^*_0(TM))$.   Since $s_{\mu}$ is supported in  $\Gamma^{\ep/2}(\gamma^*_0(TM))$, it extends to a section of $\Loopv$.   Consider a cycle $\sigma:P\rightarrow M$. We define  $$\Sigma : P\times [0,1] \rightarrow \Loopv$$ mapping $(p,t)$ to $\sigma(p)+ts_{\mu}(\pi(\sigma(p)))$. 
\begin{lemma}
$\Sigma$ defines a cycle with boundary.
\end{lemma}
\begin{proof}
We need to check all the axioms.  We can focus on points in $P$ that map to the local chart  $\Gamma^{\ep}(\gamma^*_0(TM))\oplus \Gamma_{L^2}(\gamma^*_0(TM))$ since this chart contains the support of the section $s_\mu$.  We let $\sigma(p)=(\gamma_p, v_p)$.
We have a bound on the $L^2_1$ norm of $\gamma_p$ inside the chart.  This gives a bound on $E(\gamma_p)$.    Also, note that  $$\CSD(\Sigma(p,t))=E(\gamma_p)-\frac{1}{2}|v_p+ts_\mu(p)|^2=E(\gamma_p)-\frac{1}{2}|v_p|^2-\frac{1}{2}|ts_\mu(p)|^2-\langle v_p, ts_\mu(p) \rangle$$ 
Axiom 1:  The bound on $L^2_1$ of $\gamma_p$ in the chart implies the bound on $E(\gamma_p)$ which in turn bounds $\CSD$ above.  Suppose that $(\gamma_{p_i},v_{p_i}+t_is_\mu(p_i))$ converges weakly to $(a,b)$.  We need to show that $$\CSD(a,b)\geq \liminf \CSD(\gamma_{p_i},v_{p_i}+t_is_\mu(p_i))$$  By passing to a subsequence, we can assume that $\lim (|\eta_{p_i}|_{L^2_1}^2)=\tau_\infty$ and $\lim t_i=t_\infty$.  Therefore, $t_is_\mu(p_i)=t_i\rho(|\eta_{p_i}|_{L^2_1}^2)\mu$ converges strongly $t_\infty \rho(\tau_\infty)\mu$. Let $v_\infty$ be the weak limit of $v_{p_i}$. The strong convergence of $t_i s_\mu(p_i)$ implies that $$\lim( -\frac{1}{2}|t_i s_\mu(p_i)|^2-\langle v_{p_i}, t_i s_\mu(p_i) \rangle)=-\frac{1}{2}|t_\infty \rho(\tau_\infty)\mu|^2-\langle v_\infty, t_\infty \rho(\tau_\infty)\mu \rangle$$ Since $E(a)-\frac{1}{2}|v_\infty|^2 \geq \liminf(E(\gamma_{p_i})-\frac{1}{2}|v_{p_i}|^2)$, we have $\CSD(a,b)\geq \liminf \CSD(\gamma_{p_i},v_{p_i}+t_i s_\mu(p_i))$ as desired.\\\\
Axiom 2:  Inside the chart, we have a uniform $L^2_1$ bound on $\gamma_{p_i}$.  Passing to a subsequence we can assume $\gamma_{p_i}$ converges weakly.  If we assume that $\CSD(\Sigma(p_i,t_i))$ is bounded below,  we have that $|v_{p_i}+t_is_\mu(p_i)|^2$ is bounded.  This implies that $v_{p_i}$ has bounded $L^2$ norm.  By passing to a subsequence, we can assume that $t_i s_\mu(p_i)$ is strongly converging in $L^2$.  Therefore, $(\gamma_{p_i},v_{p_i}+t_i s_\mu(p_i))$ has a weakly convergent subsequence.\\\\
Axiom 3: By our proof of Axiom 1, we see that  $\CSD(a,b)\geq \lim \CSD(\gamma_{p_i},v_{p_i}+s_\mu(p_i))$ implies that $\CSD(a,v_\infty)=\lim (\CSD(\gamma_{p+i},v_{p_i}))$.  Therefore, by passing to a subsequence, we can assume that $p_i$ converges in $P$.  By the compactness of $[0,1]$, we can assume that $t_i$ converges as well.\\\\
Axiom 4.  In the local chart, $im(\pi^+\circ D\Sigma) = im(\pi^+\circ D\sigma)$.  Therefore, it is compact.  To show that $\pi^-\circ D\Sigma$ is Fredholm we can restrict to a slice $t=t_0$  In this case, at a point $(p,t_0)$ we have $$\pi^-\circ D\Sigma_{|t_0}=\pi^-\circ (D\sigma+t_0\mu L)$$  Here $L:TP\rightarrow \Ar$ is the linear map induced by $p \mapsto \rho(|\eta_{p}|_{L^2_1}^2)$.  Since the image of $t_0\mu L$ is 1-dimensional the result follows.
\end{proof}
Using this local perturbation, we can now show that we can perturb a cycle to be transverse to $\Loop$.  Consider a cycle $\sigma:P\rightarrow \Loopv$.  Since $P\cap \Loop$ is compact, we can find finitely many sections $s_{\mu_1}, \cdots s_{\mu_k}$ such that $\Sigma: P\times [0,1]^k\rightarrow \Loopv$ with $$\Sigma(p,t_1,\cdots t_k)=\sigma(p)+t_1s_{\mu_1}(p) \cdots t_k s_{\mu_k}(p)$$ is transverse to $\Loopv$ at $P\times 0$.  We claim that be restricting $\Sigma$ to $[0,\ep]^k$, we can assume that $\Sigma$ is transverse to $\Loop$.  By contradiction, assume that there exists $p_i$ and $t^i_j$ with $t^i_j\rightarrow 0$ such that $\Sigma$ is not transverse to $\Loop$ at $p_i\times t^i_1\times \cdots t^i_k$.  By the previous lemma $\Sigma$ defines a cycle with corners,  therefore $P\times [0,1]^k\cap \Loop$ is compact.  Thus, we can assume that $p_i$ converge to $p_\infty$.  However, $\Sigma$ is transverse to $\Loop$ at $p_\infty \times 0$.  Therefore, $\Sigma$ is transverse in a small neighborhood of  $p_\infty \times 0$.  This contradicts the choice of $t^i_j$.  We have:
\begin{theorem}
Given a cycle $\sigma:P\rightarrow \Loopv$ (possibly with boundary), we can choose sections $s_{\mu_j}$ such that $\Sigma:P\times [0,\ep]^k$ is transverse to $\Loop$.  Therefore, a generic choice of $\{t_j\} \in [0,\ep]^k$ implies that $\sigma'=\Sigma_{|P\times \{ t_j\}}$ is transverse to $\Loop$.  Furthermore, assume that $P$ has no boundary.  Then, for any two such choices of $\sigma'$, $\sigma''$, there exists a cycle with boundary $H:P\times [0,1] \rightarrow \Loopv$ transverse to $\Loop$ with $\partial H=\sigma'\sqcup -\sigma''$.
\end{theorem}                
\begin{proof}
This is a standard application of Sard's theorem.  Namely, $f:P\times [0,\ep]^k \cap \Loop \rightarrow [0,\ep]^k$ is a smooth map of finite dimensional manifolds.  We can take $\{ t_j\}$ to be a regular value of this map.  $\Sigma_{|P\times \{ t_j\}}$ gives the desired $\sigma'$.  Given any two such choices, we can connect them by a path $g:[0,1]\rightarrow  [0,\ep]^k$ such that $g$ is transverse to $f$.  The map $H:P\times [0,1]\rightarrow \Loopv$ with $H(p,t)=\Sigma(p,g(t))$ provides the desired homotopy.      
\end{proof}
With this transversality result in place, we can now define a map $$i^*:\Bord_*(\Loopv) \rightarrow \Omega_*(\Loop)$$ by intersecting transverse representatives of the bordism class with $\Loop$. 
\begin{lemma}
$i^*$ is well defined.
\end{lemma}                              
\begin{proof}
By our perturbation results, we can always find a represenative $\sigma:P\rightarrow \Loopv$ in any bordism class that is transverse to $\Loop$.  This does not depend on the choice of $\sigma$ since any such choices are bordant by a transverse bordism.  Finally, if $\sigma$ is trivial, any perturbation of it will be trivial as well.  For if $\sigma\cong \partial \tau$, we can find $\{t_i\}$ as above such that but $\sigma'$ and $\tau'$ are transverse.  Since  $\partial \tau' \cong \sigma'$ we have that $\partial(\tau'\cap \Loop)\cong \sigma' \cap \Loop$  as well.  
\end{proof}
\section{Main Isomorphism}   
In the last section we defined $i^*:\Bord_*(\Loopv) \rightarrow \Omega_*(\Loop)$ by taking generic intersections with the inclusion $i:\Loop\rightarrow \Loopv$.  Let $\pi:\Loopv \rightarrow \Loop$ be the projection.  Given $\sigma:P\rightarrow \Loop$, with $\sigma\in \Omega_*(\Loop)$, let $$\pi^*(\sigma):\pi^{-1}P\rightarrow \Loopv$$ be the pullback.  
\begin{lemma}
$\pi^*:\Omega_*(\Loop)\rightarrow \Bord_*(\Loopv)$.  
\end{lemma}
\begin{proof}
$\pi^{-1}P$ consists of pairs $(p,v)$ with $v\in \Gamma(\sigma(p)^*(TM))$.  Note that $\CSD(\sigma(p),v)=E(\sigma(p))-\frac{1}{2}|v|^2$.\\\\
Axiom 1: Since $P$ is compact, we have a bound on $E(\sigma(p))$.  Therefore, $\CSD$ is bounded above on the image of $i^*(\sigma)$.  Since $-|v_i|^2$ can only rise in a weak limit we have the desired semi-continuity.\\\\  
Axiom 2: A bound on $\CSD(\sigma(p_i),v_i)$ gives us a bound on $|v_i|^2$ since $E(\sigma(p))$ is universally bounded.  Since $P$ is compact, we can assume that $p_i$ is convergent.  The bound on the norm of $v_i$ implies that $(p_i,v_i)$ has a weakly convergent subsequence.  \\\\
Axiom 3: Assume $(\sigma(p_i),v_i)$ converges weakly to $(a,b)$ with $\lim \CSD(\sigma(p_i),v_i)=\CSD(a,b)$. By the proof of Axiom 1, we must have $\lim \CSD(\sigma(p_i))=\CSD(a)$ We may assume that $p_i$ is strongly convergent.  This implies that $\lim |v_i|^2=|b|^2$.  Therefore, $v_i$ converges strongly as well.\\\\
Axiom 4:  Take a chart $\Gamma^\ep(\gamma^*(TM))\oplus \Gamma_{L^2}(\gamma^*(TM))$.  Locally, we may describe $i^{-1}(P)$ as $P\times \Gamma_{L^2}(\gamma^*(TM))$ mapping $(p,v)$ to $(\sigma(p),v)$.  Since $P$ is compact, it is clear that $\pi^+\circ Di^{*}(\sigma)=D\sigma $ is compact while $\pi^-\circ Di^{*}(\sigma)$ is an epimorphism with kernel $T_p(P)$.  
\end{proof}
By definition, we have $i^*(\pi^*(\sigma))=\sigma$.  Our goal now is to define a homotopy between $\pi^*\circ i^*$ and $Id$ on $\Bord_*(\Loopv)$.  Take $\sigma:P\rightarrow \Loopv$ to be a cycle transverse to $\Loop$.      
\begin{definition}
Let $H(P)\subset P\times [0,1]\times \Loopv$  be of triples $(p,t,(\gamma, v))$ with $\sigma(p)=(\gamma, tv)$. We let $H(\sigma):H(P)\rightarrow \Loopv$ be the projection to the last factor.  This is a smooth Hilbert manifold since $\sigma$ is assumed to be transverse to the zero section.  Note that $\partial H(\sigma)=\sigma\sqcup -\pi^*\circ i^*(\sigma)$.    
\end{definition}
\begin{lemma}
$H(\sigma)$ defines a bordism in $\Bord_*(\Loopv)$ between $\sigma$ and $\pi^*\circ i^*(\sigma)$. 
\end{lemma}
\begin{proof}
We need to verify all the axioms of a cycle. \\\\
Axiom 1: Given a triple $(p,t,(\gamma, v))\in H(P)$, note that $$\CSD(\sigma(p))=E(\gamma)-\frac{1}{2}|tv|^2\geq E(\gamma)-\frac{1}{2}|v|^2 =\CSD(\gamma,v)$$  Therefore, $\CSD$ is bounded above on the image of $H(\sigma)$.  Now, given a sequence $(p_i,t_i,(\gamma_i,v_i))$ assume that $(\gamma_i,v_i)$ is weakly convergent to $(\gamma_\infty,v_\infty)$.  By passing to a subsequence, we can assume that $t_i$ is converging.  We have $\sigma(p_i)=(\gamma_i, t_iv_i)$ weakly converging as well.  Since $\sigma$ is a cycle, we have $$\liminf \CSD(\gamma_i,t_iv_i)\leq \CSD(\gamma_\infty,t_\infty v_\infty)$$  By passing to a subsequence, we may assume that  $E(\gamma_i)$ and $|v_i|^2$ converge.   We have $$\CSD(\gamma_\infty, v_\infty )-\CSD(\gamma_\infty, t_\infty v_\infty )=-\frac{1-t_\infty}{2}|v_\infty|^2$$  and $$\lim \CSD(\gamma_{i}, v_i )-\lim \CSD(\gamma_i, t_i v_i )=-\frac{1-t_\infty}{2}\lim |v_i|^2$$ Since $\lim |v_i|^2\geq |v_\infty|^2$ we get $$\CSD(\gamma_\infty,v_\infty)- \lim \CSD(\gamma_i,v_i)=\CSD(\gamma_\infty,t_\infty v_\infty)- \lim \CSD(\gamma_i,t_i v_i)+\frac{1-t_\infty}{2}(\lim |v_i|^2-|v_\infty|^2)\geq 0$$ 
Axiom 2:  Since $\CSD(\sigma(p_i))\geq \CSD(\gamma_i,v_i)$, a lower bound on $\CSD(\gamma_i,v_i)$ implies that $\sigma(p_i)$ has a weakly convergent subsequence.  Therefore, $\gamma_i$ has bounded energy.  It follows that $\CSD(\gamma_i,v_i)$ gives a bound on the norm of $v_i$.  Therefore, $(\gamma,v_i)$ has a weakly convergent subsequence.\\\\
Axiom 3: We assume that $(\gamma_i,v_i)$ is weakly convergent to $(\gamma_\infty,v_\infty)$ and that $$\lim \CSD(\gamma_i,v_i)=\CSD(\gamma_\infty,v_\infty)$$  By passing to a subsequence, we may assume that $t_i\rightarrow t_\infty$.  By the proof of Axiom 1,  we have $$\lim \CSD(\gamma_i,v_i)-\CSD(\gamma_\infty,v_\infty)\geq \CSD(\gamma_\infty,t_\infty v_\infty)-\lim \CSD(p_i,t_iv_i)\geq 0$$  Therefore, $\CSD(\gamma_\infty,t_\infty v_\infty)=\lim \CSD(p_i,t_iv_i)$.  Since $\sigma$ is a cycle, we can pass to a subsequence where $p_i$ converge in $P$.  This implies that $\gamma_i$ converges strongly and hence, $|v_i|$ converges to $|v_\infty|$.  Therefore, $v_i$ converges strongly as well. \\\\
Axiom 4:   We may work in a fixed $t_0$ slice since this will only change the kernel/cokernel by a 1-dimensional subspace.  For $t_0>0$, $H(\sigma)$ restricted to $t_0$ is simply the map $t_0^{-1}\sigma$.  Since multiplication by $t_0^{-1}$ preserves the polarization, we have that $\pi^+(t^{-1}_0\circ \sigma )$ is compact while   $\pi^-(t^{-1}_0\circ \sigma )$ is Fredholm.  For $t_0=0$, we may identify the $t_0$ slice with $\pi^*(i^*(\sigma))$.  The previous lemma implies the conclusion for $t_0=0$.
\end{proof}
Finally, since any element of $\Bord_*(\Loopv)$ has a representative that is transverse to the zero section, we have that $\pi^*\circ i^*(\sigma)=Id$.  We summarize our result as:
\begin{theorem}
$\pi^*:\Omega_*(\Loop)\rightarrow \Bord_*(\Loopv)$ is a grading preserving isomorphism.
\end{theorem}

\section{Viterbo's Theorem and the Thom Isomorphism}
The last section established a direct isomorphism between the bordism groups of finite dimensional manifolds and that of certain semi-infinite cycles.  It is possible to establish a similar isomorphism between the cobordism versions of these theories.  In fact, the isomorphism in this case is even more direct as it is induced by the identity map.  Let $$i_*:\Omega^*(\Loop)\rightarrow \Bord^*(\Loop)$$ be the map sending $\tau:Q\rightarrow \Loop$ to $i\circ \tau:Q\rightarrow \Loopv$.  Let $$\pi_*:\Bord^*(\Loopv)\rightarrow \Omega^*(\Loop)$$ be the map sending $\tau:Q\rightarrow \Loopv$ to $\pi\circ \tau:Q \rightarrow \Loop$. We have:
\begin{lemma}
The map $i_*$ is well defined.  
\end{lemma}
\begin{proof}
Since $\CSD_{|\Loop}=E$, the $\CSD$ is bounded below by 0.  Suppose, $\tau(q_i)$ converges weakly.  $\CSD$ can only drop in weak limit since its the square of a norm on $\Loop$.  If $\CSD$ does not drop, we must have $L^2_1$ convergence of $\tau(q_i)$ which by properness of $\tau$ implies the convergence of $q_i\in Q$ (up to a subsequence).  Finally, we note that a bound on $\CSD$ implies a bound on $E$ which in turn implies the weak compactness of $\tau(q_i)$.   
\end{proof}
We will check that $\pi^*$ is well defined below.  Note that on the level of sets, $\pi_*\circ i_*=Id$.  To check that $i_*\circ \pi_*=Id$ we need to perform a homotopy:
\begin{definition}
Given $\tau:Q\rightarrow \Loopv$ in $\Bord^*(\Loopv)$, let $H(\tau):Q\times [0,1]\rightarrow \Loopv$ be the map taking $(q,t)$ to $t\cdot \tau(q)$.  
\end{definition} 

\begin{lemma}
$H(\tau)$ is a cobordism between $\tau$ and $i_*(\pi_*(\tau))$.  
\end{lemma}
\begin{proof}
As before, we need to check that all the axioms are met. \\\\
Axiom 1: We have $$\CSD(H(q,t))=E(\gamma_q)-\frac{1}{2}t|v_q|^2$$  Therefore, $\CSD(H(q,t))\geq \CSD(\tau(q))$.  This implies that $\CSD$ is bounded below on the image of $H$. Let $\tau(q_i)=(\gamma_i,v_i)$. Suppose that $(\gamma_{i},t_iv_{i})$ converges weakly.  This gives a upper bound on $\CSD(H(q_i,t_i))$ and hence $\CSD(\tau(q))$. Since, $E(\gamma_{i})$ is bounded, this bound implies that $|v_{i}|^2$ is uniformly bounded.  By passing to a subsequence, we can assume that $v_{i}$ weakly converges to $v_\infty$ and $\gamma_{i}$ weakly converges to $\gamma_\infty$.  We can also assume $t_i\rightarrow t_\infty$, and that the limits of $E(\gamma_{i})$ and $|v_i|^2$ exist.  Since $\tau$ is a cycle, we have $\CSD(\gamma_\infty,v_\infty )\leq \lim \CSD(\gamma_{i},v_{i})$.  In addition, $$\CSD(\gamma_\infty,t_\infty v_\infty)- \CSD(\gamma_\infty,v_\infty)=\frac{1-t_\infty}{2}|v_\infty|^2$$ and $$\lim \CSD(\gamma_{i},t_i v_i)- \lim \CSD(\gamma_i,v_i)=\frac{1-t_\infty}{2}\lim |v_i|^2$$  Since $\lim |v_i|^2 \geq |v_\infty|^2$, we have        $\lim \CSD(\gamma_{i},t_i v_i)\geq \CSD(\gamma_\infty,t_\infty v_\infty)$ as desired.  \\\\
Axiom 2: Suppose we have an upper bound on $\CSD(H(q_i,t_i))$.  This implies an upper bound on $\CSD(\tau(q_i))$.  Therefore, since $\tau$ is a cycle,  $(\gamma_{i},v_{i})$ has a weakly convergent subsequence.  We can assume that $t_i$ is convergent.  Therefore, we get a weakly convergent sequence $(\gamma_{i},t_iv_{i})$.  \\\\
Axiom 3:  By the proof of Axiom 1, the drop  $\lim \CSD(\gamma_{i},v_{i})-\CSD(\gamma_\infty,v_\infty)$ is a most equal to the drop  $\lim \CSD(\gamma_{i},t_iv_{i})-\CSD(\gamma_\infty,t_\infty v_\infty)$.   Therefore, if we assume that   $\lim \CSD(\gamma_{i},t_iv_{i})=\CSD(\gamma_\infty,t_\infty v_\infty)$, we get $\lim \CSD(\gamma_{i},v_{i})=\CSD(\gamma_\infty,v_\infty)$.  After passing to a subsequence, we obtain strong convergence of $q_i$ and $t_i$.  
     
\end{proof}
A proof of the (co)homology version of Viterbo's theorem can be carried out by virtually identical arguments once the relevant geometric homology theory is setup.  Namely, one considers mappings $P\rightarrow \Loopv$ with $P$ a Hilbert manifold with corners satisfying the same topological assumptions.  One can define a geometric version of singular semi-infinite homology based on such cycles \cite{MorseHomology}.  Let use denote the resulting groups by $HF_*(\Loopv)$ and the corresponding cohomology groups by $H^*(\Loopv)$.  See \cite{GeometricHomology} for a finite dimensional analogue of this construction.  The same topological arguments as in this paper give isomorphisms $$HF_*(\Loopv)\cong H_*(\Loop); HF^*(\Loopv)\cong H^*(\Loop)$$  
Another issue not addressed in this paper is the connection of this semi-infinite version of Floer theory to the traditional Morse-theoretic approach.  To connect semi-infinite homology to Morse-Floer theory  one needs to show that the $L^2$ gradient flow used in the Morse-theoretic construction gives rise to a well-defined map on the geometric chains.  This is addressed in \cite{Gradient}.  With these prerequisites, a semi-infinite version of the Morse Homology theorem in \cite{MorseHomology} proves that if the action functional is Morse-Smale, the Morse-Floer homology is isomorphic to the geometric homology of semi-infinite cycles. \\\\
Finally, let us mention that the cohomology version of Viterbo's theorem is closely related to a geometric version of the Thom isomorphism for vector bundles.  Let us pause to briefly review the geometric version of this statement.  Given an oriented finite dimensional vector bundle $$\pi: V\rightarrow M$$ of dimension $n$, the Thom isomorphism can be stated as follows.  First, one proves the existence of a Thom class $U\in H^n_{c}(V)$.  Here $ H^n_{c}(V)$ denotes cohomology with compact supports.  Given a class $v\in H^*(M)$, we may form $T(v)=U\cup \pi^*(v)\in H^*_c(V)$.  Thom's theorem is the claim that $T$ is an isomorphism.  \\\\
Let us recast this in geometric terms.  For details, consult \cite{GeometricHomology}.  Cohomology classes in $M$ are represented by maps $\tau: Q\rightarrow M$ where $Q$ is a smooth compact manifold with corners.  To get an oriented class, we must further require $\tau$ to be an oriented map.  $\pi^*(\tau)$ corresponds to the pullback $$\tau^*(V)\rightarrow V$$  The Thom class $U$ is geometrically represented by the zero section $i:M\rightarrow V$.  One way to see this is to observe that the pullback of $M$ to any fiber $V_p$ is just the origin.  This is the geometric representative of the generator of $H^n_c(V_p)$.  Now, intersection is the geometric analogue of cup products.  Therefore, $$T(\tau)=\tau^*(V)\cap M=i \circ \tau:P\rightarrow V$$
Thus, we see that geometrically, the map $T$ simply corresponds to the inclusion $T(\sigma)=i\circ \sigma$.\\\\
In the case of Viterbo's theorem, we have an infinite dimensional Hilbert bundle $$\pi:\Loopv \rightarrow \Loop$$  Just like in Thom's theorem, the map in the isomorphism $$i_*:\Omega^*(\Loop)\rightarrow \Bord^*(\Loopv)$$ is also induced by inclusion.  Therefore, we may view the cobordism version of Viterbo's theorem (or more precisely its cohomological counterpart) as a semi-infinite analogue of the Thom isomorphism theorem.  \\\\

\newpage
\bibliography{thesisbib}
\bibliographystyle{plain}

\end{document}